\documentclass{amsproc}

\newtheorem{Theorem}{Theorem}

\newtheorem{theorem}{Theorem}[section]
\newtheorem{lemma}[theorem]{Lemma}
\newtheorem{proposition}[theorem]{Proposition}
\newtheorem{corollary}[theorem]{Corollary}
\theoremstyle{definition}

\newtheorem{example}[theorem]{Example}

\newtheorem{assumption}[theorem]{Assumption}

\theoremstyle{remark}
\newtheorem{remark}[theorem]{Remark}

\numberwithin{equation}{section}

\newcommand{\abs}[1]{\lvert#1\rvert}

\newcommand{\R}{{\mathbb R}}

\newcommand{\N}{{\mathbb N}}

\newcommand{\RR}{{\mathcal R}}

\begin{document}

\title[Standard Reeb flow]
{A characterization of the standard Reeb flow}


\author{Shigenori Matsumoto}
\address{Department of Mathematics, College of
Science and Technology, Nihon University, 1-8-14 Kanda, Surugadai,
Chiyoda-ku, Tokyo, 101-8308 Japan
}
%\curraddr{Department of Mathematics, College of
%Science and Technology, Nihon University, 1-8-14 Kanda, Surugadai,
%Chiyoda-ku, Tokyo, 101-8308 Japan}
\email{matsumo@math.cst.nihon-u.ac.jp
}
\thanks{The author is partially supported by Grant-in-Aid for
Scientific Research (C) No.\ 20540096.}
%    General info


\date{\today }

\begin{abstract}
Among the topological conjugacy classes of the continuous flows
$\{\phi^t\}$ whose orbit foliations are the planar Reeb foliation, there
is one special class called the standard Reeb flow. We show that
$\{\phi^t\}$ is conjugate to the standard Reeb flow if and only if
$\{\phi^t\}$ is conjugate to $\{\phi^{\lambda t}\}$ for any $\lambda>0$.
\end{abstract}

\maketitle

\section{Introduction}
Let
$$
P=\{(\xi,\eta)\mid \xi\geq0,\eta\geq0\}-\{(0,0)\}.$$
A nonsingular flow $\{\Phi^t\}$ on $P$ defined by
$$
\Phi^t(\xi,\eta)=(e^t\xi,e^{-t}\eta)$$
is called the {\em standard Reeb flow}. In this note the oriented
foliation $\RR$ whose leaves are the orbits of $\{\Phi^t\}$
with the orientation given by the time direction
 is called the {\em Reeb foliation}. A continuous flow on $P$
with orbit foliation
$\RR$ is called an $\RR$-flow.
The topological conjugacy classes of $\RR$-flows $\{\phi^t\}$
are classified in \cite{L}
in the following way. Let 
$\gamma_1:[0,\infty)\to P$ (resp.\ $\gamma_2:[0,\infty)\to P$)
be a continuous path such that $\gamma_1(0)\in\{\xi=0\}$
(resp.\ $\gamma_2(0)\in\{\eta=0\}$) which intersects every interior leaf of $\RR$
at exactly one point.
Then one can define a continuous function 
$$f_{\{\phi^t\},\gamma_1,\gamma_2}:(0,\infty)\to \R$$
by setting that $f_{\{\phi^t\},\gamma_1,\gamma_2}(x)$
is the time needed for the flow $\{\phi^t\}$ to
move from the point $\gamma_1(x)$ until it 
reaches a point on the curve $\gamma_2$.
Then $f_{\{\phi^t\},\gamma_1,\gamma_2}$ belongs to the
following space 
$$E=\{f:(0,\infty)\to\R\mid f\ \ \mbox{is continuous and}
\ \ \lim_{x\to 0}f(x)=\infty\}.$$
Of course $f_{\{\phi^t\},\gamma_1,\gamma_2}$ depends upon the
choices of $\gamma_1$ and $\gamma_2$. There are two umbiguities,
one coming from the parametrization of $\gamma_1$, and the other
coming from the positions of $\gamma_1$ and $\gamma_2$.
Let $H$ be the space of homemorphisms of $[0,\infty)$ and $C$
the space of continuous functions on $[0,\infty)$.
Define an equivalence relation $\sim$ on $E$ by
$$
f\sim f'\Longleftrightarrow f'=f\circ h+k,\ \  \exists h\in H,\
\ \exists k\in C.$$
Then clearly the equivalence class of $f_{\{\phi^t\},\gamma_1,\gamma_2}$
does not depend on the choice of $\gamma_1$ and $\gamma_2$. Moreover
it is an invariant of the topological conjugacy classes of
$\RR$-flows. Thus if we denote
by $\mathcal E$ the set of the topological conjugacy classes of 
the $\RR$-flows, then there is a well defined map
$$
\iota:\mathcal E\to E/\sim.$$
The main result of \cite{L} states that
$\iota$ is a bijection. In particular any $f\in E$ is obtained
as $f=f_{\{\phi^t\},\gamma_1,\gamma_2}$ for some $\RR$-flow $\{\phi^t\}$
and paths $\gamma_i$.

Clearly any strictly monotone function of $E$ belongs to a single
equivalence class, and this corresponds to the standard Reeb 
flow $\{\Phi^t\}$. The purpose of this note is to show the following
characterization of the standard Reeb flow.

\begin{Theorem} \label{t}
An $\RR$-flow $\{\phi^t\}$ is topologically conjugate to
the standard Reeb flow $\{\Phi^t\}$ if and only if $\{\phi^{\lambda t}\}$
is topologically conjugate to $\{\phi^t\}$ for any $\lambda>0$.
\end{Theorem}
 
Of course the only if part is immediate. We shall show the if part
in the next section.

\begin{remark} \label{r10}
A single $\lambda$ is not enough for Theorem \ref{t}.
In fact there is an $\RR$-flow $\{\phi^t\}$ not
topologically conjugate to $\{\Phi^t\}$ such that
$\{\phi^{2t}\}$ is topologically conjugate to $\{\phi^t\}$.
This will be given in Example \ref{example} below.
\end{remark}

The author wishes to express his hearty thanks to the anonymous
referees whose valuable comments are indeed helpful for the improvement
of the paper

\section{Proof of the if part}
The equivalence class of $f\in E$ is determined by how 
$f(x)$ oscilates while it tends to $\infty$ as $x\to 0$.
So to measure the degree of oscilation of
$f\in E$, define a nonnegative valued continuous function
$f^*$ defined on $(0,1]$ by
$$
f^*(x)=\max(f\vert_{[x,1]})-f(x).
$$
Then we have the following lemma.

\begin{lemma} \label{easy}
(1) If $\lambda>0$, then $(\lambda f)^*=\lambda f^*$.

(2) If $c$ is a constant, then $(f+c)^*=f^*$.

(3) If $h\in H$, then there is $0<a<1$ such that
$(f\circ h)^*=f^*\circ h$ on $(0,a)$.

(4) If $k\in C$ and $x\to 0$, then $(f+k)^*(x)-f^*(x)\to 0$.

(5) There is a sequence $\{x_n\}$ tending to $0$
such that  $f^*(x_n)=0.$
\end{lemma}

{\sc Proof}. Points (1) and (2) are immediate. 
To show (3) notice that
\begin{eqnarray*}
&(f\circ h)^*(x)=\max(f\vert_{[h(x),h(1)]})-f(h(x))\ \mbox{ and }\\
&f^*\circ h(x)=\max(f\vert_{[h(x),1]})-f(h(x)).
\end{eqnarray*}
Since $f(x)\to\infty$ ($x\to0$),  both maxima coincide for small $x$.

Let us show
(4). By (2) we only need to show (4) assuming that
$k(0)=0$. Now given $\epsilon>0$, there is $\delta>0$
such that if $0<x<\delta$, then $\abs{k(x)}<\epsilon$.
Choose $\eta>0$ small enough so that if $0<x<\eta$,
then we have
$$
f(x)\geq\max(f\vert_{[\delta,1]})\ \ \mbox{and}\ \
(f+k)(x)\geq\max((f+k)\vert_{[\delta,1]}).
$$
This implies that for $x\in(0,\eta)$,
$$
\abs{f^*(x)-(f+k)^*(x)}
\leq\abs{f(x)-(f+k)(x)}+\abs{\max((f+k)\vert_{[x,\delta]})
-\max(f\vert_{[x,\delta]})}<2\epsilon.
$$
This shows (4). Finally (5) follows from the assumption $f(x)\to\infty$
as $x\to0$.
\qed

\bigskip

For $f\in E$ define an invariant $\sigma(f)=\limsup_{x\to 0}f^*(x)$
which takes value in $[0,\infty]$. In fact $\sigma(f)$ coincides
with the invariant $\mathcal A(f)$ defined in \cite{L} and used
to show that  $\mathcal E$ is uncountable.

\begin{lemma}\label{sigma}
Assume $f,f'\in E$ and $\lambda>0$.

(1) We have $\sigma(\lambda f)=\lambda\sigma(f)$.

(2) If $f\sim f'$, then $\sigma(f)=\sigma(f')$. 
In particular $f$ corresponds to the standard Reeb flow 
if and only if $\sigma(f)=0$.
\end{lemma}

{\sc Proof}. Clearly (1) follows from Lemma \ref{easy} (1),
while the first statement of (2) is an easy consequence of Lemma
\ref{easy}
(3) and (4). 
To show the last statement, assume $\sigma(f)=0.$ Extend the function
$f^*$ defined on $(0,1]$ to $[0,\infty)$ by letting
$$f^*=0\ \mbox{ on }\ \{0\}\cup(1,\infty).$$ 
Since $\sigma(f)=0$, $f^*$ is continuous, i.\ e.\ $f^*\in C$. 
Thus $f\sim f+f^*$, and the latter is (weakly) monotone near $0$.
Still adding a suitable function, $f$ is seen to be equivalent to
$g$ which is strictly monotone on the whole $(0,\infty)$ such that
$g(x)\to0$ ($x\to\infty$). 
Clearly such functions are mutully equivalent by
a pre-composition of some $h\in H$, and correspond to the standard
Reeb flow $\{\Phi^t\}$.
\qed

\bigskip

Now since 
\begin{equation} \label{e101}
f_{\{\phi^{\lambda t}\},\gamma_1,\gamma_2}=
\lambda^{-1} f_{\{\phi^t\},\gamma_1,\gamma_2},
\end{equation}
for $\lambda>0$,
Theorem \ref{t} reduces to the following proposition.

\begin{proposition} \label{p}
If $f\in E$ and $f\sim\lambda f$ for any $\lambda>0$, then $\sigma(f)=0$.
\end{proposition}

The rest of the paper is devoted to the proof of Proposition \ref{p}.
But before starting, let us mention an example
for Remark \ref{r10}.

\begin{example} \label{example}
By (\ref{e101}) and the main result of \cite{L},
it suffices to construct a function $f\in E$ such that $f(x/2)=2f(x)$ and
that $\sigma(f)=\infty$. Set for example
$$
f(x)=\frac{1}{x}2^{\sin(2\pi\log_2x)}.$$
\end{example}

The following lemma, roughly the same thing as
the linearization in one dimensional local
dynamics, plays a crucial role in what follows.

\begin{lemma} \label{key}
Assume $f\in E$ satisfies $\lambda f=f\circ h+k$ for some $h\in H$,
$k\in C$ and $\lambda>1$. Then 
$0$ is an attracting fixed point of $h$
and there exists $f_\infty\in E$ such that
$f_\infty-f\in C$, $\lambda f_\infty=f_\infty\circ h$
and $f_\infty(x)\to 0$ ($x\to\infty$).
\end{lemma}

{\sc Proof}. 
Any equivalence class of $E$ has a representative $f$
such that
\begin{equation} \label{e1}
f\vert_{[1,\infty)}\ \ \mbox{is bounded}.
\end{equation}
So it is no loss of generality to assume that the function $f$ in
the lemma satisfies (\ref{e1}). We can also assume that $k(0)=0$,
by adding a suitable constant to $f$ if necessary.
Choose $a'\in(0,1)$ so that if $a\in(0,a')$,
$$
f(a)>\frac{2}{\lambda-1}\max(\abs{k}\vert_{[0,1]}).$$
Then we have
\begin{equation} \label{e100}
f\circ h(a)>\frac{\lambda+1}{2}f(a),\ \ \forall a\in(0,a').
\end{equation}
If $a$ is sufficiently
near $0$, we have
$$
f(a)>\sup(f\vert_{[1,\infty)}).$$ 
If furthermore $f^*(a)=0$, then 
$$
\{x\mid f(x)>f(a)\}\subset(0,a).$$
Thus (\ref{e100})
implies
 $h(a)<a$ for such $a$.
But this allows us to use (\ref{e100}) repeatedly for $h^n(a)$ 
($n=1,2,\cdots$) instead of $a$,
showing that $f\circ h^n(a)\to\infty$ as $n\to\infty$.
Clearly this implies that
$[0,a]$ is contained in the attracting domain
of an attractor $0$ of the homeomorphism $h$, showing the first point
of Lemma \ref{key}.

For the rest of the proof, let us divide the argument into
two cases according to the dynamics of $h$.
First assume that the whole line $[0,\infty)$
is the attracting domain of $0$.
Let
$$
f_n(x)=\lambda^{-n}f(h^n(x)).$$
Then we have
$$
f_{n+1}(x)-f_n(x)=-\lambda^{-n-1}k(h^n(x)),$$
showing that $f_n\to f_\infty$ uniformly on compact
subsets of $(0,\infty)$ for some
continuous function $f_\infty$.
Now since
$$\lambda f_{n+1}(x)=f_n(h(x)),$$
we have
$$
\lambda f_\infty=f_\infty\circ h.$$
We also have
$$\abs{f(x)-f_\infty(x)}\leq\sum_
{n=0}^\infty \lambda^{-n-1}\abs{k(h^n(x))}. $$
The continuity of $k$, together with the assumption $k(0)=0$,
implies that
$$\lim_{x\to0}\abs{f(x)-f_{\infty}(x)}=0,$$
showing that $f_\infty-f\in C$. 

Finally since $h^{-n}(x)\to\infty$
($n\to\infty$) and
$$f_\infty\circ h^{-n}(x)=\lambda^{-n}f_\infty(x), \ \ \forall
x\in(0,\infty),
$$
we have $f_\infty(x)\to 0$ ($x\to\infty$).

Next assume there is a fixed point $b$ of $h$ such that
$(0,b)$ is an attracting domain of $0$. 
Thus we have $h^{-n}(x)\to b$ ($n\to \infty$) for any 
$x\in (0,b)$.

The same argument
as above shows the existence of a continuous function
$f_\infty$ on $(0,b)$. 
Since
$$
f_\infty\circ h^{-n}(x)=\lambda^{-n}f_\infty(x),\ \ \forall x\in(0,b),
$$
we have
$$
\lim_{x\uparrow b}f_\infty(x)=0.$$
Now extend $f_\infty$ by setting $f_\infty=0$ on $[b,\infty)$.
\qed

\bigskip
Let us start the proof of Proposition \ref{p}.
Assume $f\in E$ satisfies $f\sim 2^{1/N}f$ for any $N\in\N$.
Applying Lemma \ref{key}, $f$ can be changed
within the equivalence class to one which satisfies the condition
of $f_\infty$ for $\lambda=2$. 
We also assume for contradiction that
$\sigma(f)>0$. 

Thus the proof of Proposition \ref{p} reduces to showing
that there is no $f\in E$ which
satisfies the following assumption.

\begin{assumption} \label{a1}
A function $f\in E$ satisfies 
\begin{eqnarray}
2f=f\circ h,\ \ \exists h\in H, \ 
& f(x)\to 0\ \ (x\to\infty),\label{first}\\
2^{1/N}f-f\circ h_N\in C\ \ & \exists h_N\in H, \label{second}\ \
\forall N\geq2\ \ \mbox{and}\\
\sigma(f)>0. & \label{third}
\end{eqnarray}

\end{assumption}

Define
$$E_0=\{f\in E\mid f(x)\to0\ \ (x\to\infty)\}.$$
Henceforth all the functions dealt with will be in $E_0$,
and the following definition is more convenient.
For $f\in E_0$ define
$$
f^\sharp(x)=\max(f\vert_{[x,\infty)})-f(x).$$
Clearly $f^\sharp$ and $f^*$ are the same near $0$ and Lemma
\ref{easy} (1), (4) and (5) hold also for $f^\sharp$, while (3)
becomes stronger. In summary we have:

\begin{lemma} \label{sharp}
Assume $f, f'\in E_0$.

(1) If $\lambda>0$, then $(\lambda f)^\sharp=\lambda f^\sharp$.

(3) If $h\in H$, then $(f\circ h)^\sharp=f^\sharp\circ h$.

(4) If $f'-f\in C$ and $x\to0$, then $f^\sharp(x)-(f')^\sharp(x)\to0$.

(5) There is a sequence $\{x_n\}$ tending to $0$ such that
 $f^\sharp(x_n)=0$.
\qed
\end{lemma}

Hereafter $f$ is always to be a function satisfying Assumption \ref{a1}.
Thus we have
\begin{equation} \label{no}
2f^\sharp=f^\sharp\circ h.
\end{equation}
Fix $N$ for a while and let $h_1=h_N^N$. 
Notice that by Lemma \ref{key} both
$h$ and $h_1$ have $0$ as their attractors
and that
$$f\circ h-f\circ h_1=2f-f\circ h_1
=\sum_{\nu=0}^{N-1}2^{\frac{N-\nu-1}{N}}(2^{\frac{1}{N}}f\circ
h_N^\nu
-f\circ h_N^{\nu+1})\in C.$$
The following is an easy corollary of Lemma \ref{sharp}.

\begin{corollary} \label{small}
We have
$$
\lim_{x\to 0}\abs{f^\sharp\circ h(x)-f^\sharp\circ h_1(x)}=0.$$
\qed
\end{corollary}

Our overall strategy is to show that $f^\sharp$ is too much
oscilating in a fundamental domain of $h$, thanks to condition
(\ref{second}).
For that purpose first of all
we have to compare the dynamics of $h$ and $h_1$ near the common
attractor $0$ and to show that they have more or less the same
fundamental domains.

\begin{lemma} \label{final}
Either there exists a sequence $\{a_n\}$ such that
$a_n\to 0$ and that $h^2(a_n)\leq h_1(a_n)\leq h(a_n)$
or there exists a sequence $\{a_n\}$ such that
$a_n\to 0$ and that $h^2_1(a_n)\leq h(a_n)\leq h_1(a_n)$.
\end{lemma}

{\sc Proof}.
If there is a sequence $\{a_n\}$ such that $a_n\to 0$
and that $h(a_{n})=h_1(a_n)$, there is nothing to prove.
So there are two cases to consider.
One is when $h_1(x)<h(x)$ for any small $x$, and the other
$h_1(x)>h(x)$. 

For the moment assume the former.
In way of contradiction assume the contrary of the assertion
of the lemma.
This is equivalent to saying that $h_1(x)< h^2(x)$
for {\em any small} $x$.
For small $x$,
let $y=y(x)\in[h_1(x),x]$ be any point which gives
$\max(f^\sharp\vert_{[h_1(x),x]})$. 
Notice that $f^\sharp(y)$ can be as large as we wish by choosing
$x$ even smaller.
Then since $f^\sharp(h^2(y))=4f^\sharp(y)>f^\sharp(y)$,
the point $h^2(y)$ is contained in 
$$[h^2\circ h_1(x),h^2(x)]- (h_1(x),x]=[h^2\circ h_1(x),h_1(x)]
\subset [h^2_1(x),h_1(x)].
$$
The last inclusion follows from the assumption for a contradiction.

Put $h^2(y)=h_1(z)$ for some $z=z(x)\in[h_1(x),x]$. Then we 
have
\begin{equation} \label{e31}
f^\sharp\circ h_1(z)=4f^\sharp(y)\geq 4f^\sharp(z)\ \ \mbox{ and }\ \ 
f^\sharp\circ h(z)=2f^\sharp(z).
\end{equation}
If we choose $x$ near enough to $0$, then 
the associated $z=z(x)$ is also near, and thus
$$\abs{2f^\sharp(z)
-f^\sharp\circ h_1(z)}
=\abs{f^\sharp\circ h(z)
-f^\sharp\circ h_1(z)}
$$ 
can be arbitrarily small by Corollary \ref{small}.
Then we have 
$$f^\sharp(z)\approx \frac{1}{2}f^\sharp\circ h_1(z)=
2f^\sharp(y)\geq 1
$$
for any such $z=z(x)$.
On the other hand $z(x)$ can be arbitrarily near to $0$, and
thus (\ref{e31}) contradicts Corollary \ref{small}.

The opposite case where $h(x)<h_1(x)$ for any small $x$
can be dealt with similarly by considering $f'\in E$,
equivalent to $f$, such that $2f'=f'\circ h_1$ and $f'(x)\to0$
($x\to\infty$), instead of $f$.
\qed

\bigskip

Now fix a large number $N$ and choose $f_1\in E_0$ such that
$$f_1-f\in C,\ \ 2^{1/N}f_1=f_1\circ h_N. 
$$
The existence of such $f_1$ is
guaranteed by Lemma \ref{key} applied to $\lambda=2^{1/N}$.
We have then 
\begin{equation}\label{yes}
2^{1/N}f_1^\sharp=f_1^\sharp\circ h_N.
\end{equation}

Together with Lemma \ref{final} which asserts that
the fundamental domain of $h_N^N$ is more or
less comparable with that of $h$, this implies that $f_1^\sharp$ is oscilating
in an extremely high frequency for $N$ big. We are going to get
a contradiction from this.

We still assume (\ref{first}) for $f$. 
According to Lemma \ref{final}, there are two cases to consider.
One is when there is a sequence $a_n\to 0$ such that
$h^2(a_n)\leq h^N_N(a_n)\leq h(a_n)$, the other being 
$h^{2N}_N(a_n)\leq h(a_n)\leq h^N_N(a_n)$.

Assume for the moment that the former holds for infinitely many $N$.
Let $x_n^1$ be the largest point such that $x^1_n\leq a_n$ and
$f_1^\sharp(x_n^1)=0$. 
Notice that by Lemma \ref{sharp} (5) and the equation (\ref{yes}), we have
\begin{equation}\label{24}
x_n^1\in(h_N(a_n),a_n].
\end{equation}
Then again by (\ref{yes})
$f_1^\sharp$ vanishes at the points $x_n^\nu=h_N^{\nu-1}(x^1_n)$
for any $1\leq \nu\leq N$. Let
$y_n^1$ be any point in $[x_n^2,x_n^1]$ at which $f_1^\sharp$
takes the maximal value
and let $y_n^\nu=h_N^{\nu-1}(y^1_n)$
for $1\leq \nu \leq N-1$.
By (\ref{24}) the order of these points are as follows.
$$
h^2(a_n)<h_N^N(a_n)\leq x^{N}_n < y^{N-1}_n<\cdots<y_n^\nu<x_n^\nu<\cdots
<y^1_n<x^1_n\leq a_n.$$
Notice that $y_n^\nu$ is a point in $[x_n^{\nu+1},x_n^\nu]$
at which $f_1^\sharp$ takes the maximal value,
and
$$f_1^\sharp(y_n^\nu)=2^{(\nu-1)/N}f_1^\sharp(y_n^1).$$
We also have
\begin{equation}\label{1/2}
f_1^\sharp(y_n^\nu)\geq \frac{1}{2}
\max(f_1^\sharp\vert_{[h_N^N(a_n),a_n]}).
\end{equation}
In fact on one hand
$$\max(f_1^\sharp\vert_{[x_n^N,a_n]})
=f_1^\sharp(y_n^{N-1})= 2^{(N-2)/N}f_1^\sharp(y_n^1)\leq
2f_1^\sharp(y_n^1).
$$
On the other hand
$$\max(f_1^\sharp\vert_{[h_N^N(a_n),x_n^N]})\leq 2^{(N-1)/N}\max
(f_1^\sharp\vert_{[x_n^2,x_n^1]}) 
\leq 2f_1^\sharp(y_n^1),$$
because
$$h_N^{-N+1}[h_N^N(a_n)), x_n^N]=[h_N(a_n),x_n^1]
\subset [x_n^2,x_n^1].
$$

Henceforth we focus our attention to the other homeomorphism
$h\in H$.
There is a sequence $\{m_n\}$ of integers such that
the points $h^{-m_n}(a_n)$ belong to a fixed fundamental domain
in the basin of $0$ for $h$. Notice that $m_n\to\infty$ since
$a_n\to0$.
Passing to a subsequence if necessary, we may assume that
$$h^{-m_n}(a_n)\to a,\ \ \
h^{-m_n}(x^\nu_n)\to x^\nu\ \mbox{and}\ h^{-m_n}(y^\nu_n)\to y^\nu,$$
for some points $a$, $x^\nu$ and $y^\nu$.
There
is an ordering
$$h^2(a)\leq x^{N} \leq y^{N-1}\leq\cdots \leq y^\nu\leq x^\nu\leq\cdots
\leq y^1\leq x^1\leq a.
$$

We shall show that $f^\sharp(x^\nu)=0$ and that $f^\sharp(y^\nu)$ is bounded
away from 0 with a bound {\em independent of $N$}. 
Since these points can be taken in the same compact interval $[h^2(a),a]$,
this will contradict the continuity of $f^\sharp$.

By Lemma \ref{sharp} (4), $f_1^\sharp(x_n^\nu)=0$
implies  $f^\sharp(x_n^\nu)\leq 1$ for any
large $n$. Therefore by(\ref{no})
$$
f^\sharp(h^{-m_n}(x_n^\nu))\leq 2^{-m_n},$$
showing that $f^\sharp(x^\nu)=0.$

On the other hand since $h_N^N(a_n)\leq h(a_n)$, we have by (\ref{1/2})
$$
f_1^\sharp(y_n^\nu)\geq \frac{1}{2}\max(f_1^\sharp\vert_{[h_N^N(a_n),a_n]})
\geq \frac{1}{2}\max(f_1^\sharp\vert_{[h(a_n),a_n]}),
$$
and therefore again by Lemma \ref{sharp} (4), for any large $n$,
$$f^\sharp(y_n^\nu)\geq\frac{1}{2}\max(f^\sharp\vert_{[h(a_n),a_n]})-1.$$
Let $M=\max(f^\sharp\vert_{[h(a),a]})$ and notice that $M>0$ since
$\sigma(f)>0$ (\ref{third}) and by (\ref{no}).

For any large $n$, the interval $h^{-m_n}[h(a_n),a_n]$
is near $[h(a),a]$, and is composed of a subinterval of $[h(a),a]$
and the iterate by $h^{\pm1}$ of the complementary subinterval, and
therefore
$$
\max(f^\sharp\vert_{h^{-m_n}[h(a_n),a_n]})\geq M/2.$$
This implies by (\ref{no})
$$
\max(f^\sharp\vert_{[h(a_n),a_n]})\geq \frac{1}{2}M2^{m_n},$$
showing that for any large $n$
$$f^\sharp(y_n^\nu)\geq\frac{1}{4}M2^{m_n}-1.$$
This concludes that
$$f^\sharp(y^\nu)\geq\frac{1}{4}M,$$
as is desired.

The opposite case where 
$h^{2N}_N(a_n)\leq h(a_n)\leq h_N^N(a_n)$ ($\exists a_n\to 0$)
holds for infinitely many $N$
can be dealt with in a similar way, although the argument is not 
completely symmetric.

\end{document}